\documentclass[12pt,a4paper]{article}
\usepackage{latexsym,amsfonts,ifthen}
\begin{document}

\newtheorem{theo}{Theorem}
\newtheorem{prop}{Proposition}
\newtheorem{lem}{Lemma}
\newtheorem{efi}{Definition}
\newtheorem{coro}{Corollary}
\newcommand{\R}{\ensuremath{\mathbb R}}
\newcommand{\qed}{\hfill $\Box$}
\newcommand{\beq}{\begin{eqnarray*}}
\newcommand{\eeq}{\end{eqnarray*}}
\newcommand{\lag}{\langle}
\newcommand{\rag}{\rangle}
\newcommand{\la}{\lambda}
\newcommand{\D}{\displaystyle}
\newboolean{integral}
\newcommand{\option}[2][{}]{\ifthenelse{\boolean{integral}}{#2}{#1}}
\newcommand{\signature}{\vskip 2truecm

	\hbox{\noindent\hskip 5truecm
	\vbox{\advance\hsize by -5truecm

	\centerline{Equipe d'Analyse et de Math\'ematiques Appliqu\'ees}

	\centerline{Universit\'e de Marne-la-Vall\'ee, 2 rue de la Butte Verte}

	\centerline{93166 Noisy-Le-Grand CEDEX. France}

	\centerline{e-mail: barthe@clipper.ens.fr}}}

	}

\setboolean{integral}{false}

\title{On a reverse form of the Brascamp-Lieb inequality}
\author{Franck Barthe \footnote{This work will form part of a doctoral thesis
	under the supervision of professors B.~Maurey and A.~Pajor. Their 
	advice and encouragement have been invaluable.}}
\maketitle

\begin{abstract} We prove  a reverse form of the multidimensional
	Brascamp-Lieb inequality. Our method also gives a new way
	to derive the Brascamp-Lieb inequality and is rather
	convenient for the study of equality cases.
\end{abstract}

{\bf Introduction}
We will work on the space $\R^n$ with its usual Euclidean structure. We will
denote by $\langle , \rangle$ the canonical scalar product.
In \cite{brasl76bcyi}, H. J. Brascamp and E. H. Lieb showed that for $m \ge n$,
$p_1,\ldots,p_m >1$ and $a_1,\ldots,a_m \in \R^n$, the norm of the multilinear
operator $\Phi$ from $L_{p_1}(\R) \times \cdots \times L_{p_m}(\R)$ into 
$\R$ defined by
\[ \Phi (f_1,\ldots,f_m)= \int _{\R^n} \prod_{i=1}^m f_i (\langle x,a_i \rangle) \, dx \]
can be computed as the supremum over centered Gaussian functions $g_1,\ldots,g_m$
of 
\[ \frac{\Phi( g_1,\ldots,g_m )}  {\prod_{i=1}^m \|g_i\|_{p_i}} \cdot \]
In other words, $\Phi$ is saturated by Gaussian functions. This theorem is a very
convenient tool to derive sharp inequalities. Brascamp and Lieb applied it successfully
to prove the optimal version of Young's convolution inequality (also derived
independently and simultaneously by Beckner \cite{beck75ifa}), to rederive Nelson's
hypercontractivity. Their proof is based on a rearrangement inequality 
of Brascamp, Lieb and Luttinger \cite{brasll74grim} and on the fact that radial 
functions of a large number of variables behave like Gaussians. However, their
method left opened, except in some special cases, the multidimensional problem:
 let $m\ge n$, $p_1,\ldots,p_m >1$ and let $n_1,\ldots,n_m$ be integers; for each $i\le m$
let $B_i$ be a linear mapping from $\R^n$ into $\R^{n_i}$. Is the multilinear operator
on $L_{p_1}(\R^{n_1}) \times \cdots \times L_{p_m}(\R^{n_m})$ defined by
\[ \Psi(f_i,\ldots,f_m)=\int_{\R^n} \prod_{i=1}^m f_i (B_i x) \, dx \]
saturated by Gaussian functions?

This question was solved positively by Lieb in his article ``Gaussian kernels have 
only Gaussian maximizers'' \cite{lieb90gkgm}. The point is that $\Psi$ can be wieved
as a limit case of multilinear operators with Gaussian kernels.

In \cite{bart97iblc}, the author gave a simple proof, for functions of one real variable,
of the Brascamp-Lieb inequality and of a new family of inequalities which can be 
understood as a reverse form, or as a dual form of the Brascamp-Lieb inequalities.
These inequalities can be stated as follows: let $m \ge n$,
$p_1,\ldots,p_m >1$ and $a_1,\ldots,a_m \in \R^n$, the largest constant $E$
such that 
\[ \int\limits_{\R^n}^{\,*}   \sup_{x=\sum c_{i} \theta_i a_{i}} \prod_{i=1}^{m} f_{i}(\theta_{i}) \, dx
	\ge
   E \prod_{i=1}^m \|f_i\|_{p_i} \]
holds for all $f_i,\ldots,f_m$ is also the largest constant such that the inequality holds 
for centered Gaussian functions, where $\int^*$ is the outer integral. Again, Gaussian
 functions play an extremal role.
This new inequality was inspired by convexity theory. The strength of the Brascamp-Lieb inequality
for volume estimates of convex bodies was noticed by K. Ball (see \cite{ball89vscr}, \cite{ball91scb}
and \cite{ball91vrri}), who also remarked in \cite{ball91vrri} that a reverse inequality would
give dual results. For geometric applications of the reverse Brascamp-Lieb inequality see 
\cite{bart98} and also section III of the present paper.

\vspace{5mm}
In the first section, we prove a fully multidimensional version of the reverse
Brascamp-Lieb inequality. Our method also gives a new proof of the multidimensional
Brascamp-Lieb inequality. It is very similar to the one dimensional case and
uses a theorem of Brenier (\cite{bren87dprm}, \cite{bren91pfmr}) refined by McCann (\cite{mcca94ctig},
\cite{mcca95eumm}) on measure preserving mappings deriving from convex potentials. Notice that this
result was applied by McCann in \cite{mcca94ctig} to prove the Pr\'ekopa-Leindler inequality 
(\cite{prek73lcmf}, \cite{lein72cchi}),
which is a particular case of the reverse Brascamp-Lieb inequalities.

\vspace{5mm}
In section II, we focus on the one-dimensional case to deal in detail with equality cases.
This problem was left opened for the Brascamp-Lieb inequality because the previous
proofs were depending on limit processes. We push further the study of \cite{brasl76bcyi}
in the spirit of \cite{lieb90gkgm} to see when there is a Gaussian maximizer for the 
Brascamp-Lieb inequality (or a Gaussian minimizer for the reverse form) and whether it
is unique.

\vspace{5mm}
In section III, we study the particular case of the Brascamp-Lieb inequality which was noticed
by K. Ball \cite{ball89vscr} and which is so usefull in convexity. We state the corresponding
converse inequality. The equality cases are completely solved,
 which allows us to find new characteristic properties of simplices and parallelotopes.
 The multidimensional version of the reverse Brascamp-Lieb inequality implies a theorem similar to
Brunn-Minkowski inequality for sets that are contained in subspaces.

\vspace{5mm}
In section IV, we develop an idea of \cite{brasl76bcyi}: after proving a reverse and sharp form
 of Young's inequality, Brascamp and Lieb take limits in some parameters and rederive the 
Pr\'ekopa-Leindler inequality. We rederive Ball's version of the Brascamp-Lieb inequality and its converse 
from a generalized form of Young's inequality and its converse (which we prove by the 
method that we developed in \cite{bart98oyic}). The goal of this section is to show the unity
of the topics.

\section{Proof of the Brascamp-Lieb inequality and its converse}

We first introduce some notations. We will denote by $\mathcal S^+(\R^n)$ the set of $n\times n$
symmetric definite positive matrices. For $A \in S^+(\R^n)$ we will denote by $G_A$ the 
centered Gaussian function on $\R^n$
	$$ G_A(x)=\exp (- \langle Ax, x \rangle) \, .$$
We will also denote by $\mathcal L (\R^n,\R^m)$ the set of linear mappings from $\R^n$ to
$\R^m$, identified with $m \times n$-matrices. If $B \in  \mathcal L (\R^n,\R^m)$,
 then $B^* \in  \mathcal L (\R^m,\R^n)$ will be its Euclidean adjoint. We will work
in the set of integrable non-negative functions on $\R^n$, denoted by $L_1^+(\R^n)$.

The fully multidimensional version of the Brascamp-Lieb inequality and its converse is
	as follows:
\begin{theo}
\label{thRBL}
Let $m\ge n$ be integers. Let $(c_i)_{i=1}^m$ be positive real numbers and $(n_i)_{i=1}^m$
be integers smaller than $n$ such that
$$\sum_{i=1}^m c_i n_i =n.$$
For $i=1 \ldots n$ let $B_i$ be a linear surjective map from $\R^n$ onto $\R^{n_i}$. Assume that
$$ \cap_{i \le m} \ker \, B_i = \{ 0 \}.$$
We define two applications $I$ and $J$ on $L^+_1(\R^{n_1})\times \cdots \times L^+_1(\R^{n_m})$ as
follows: if $f_i \in  L^+_1(\R^{n_i})$, $i=1, \ldots ,m$ then

$$ I((f_i)_{i=1}^m)= \int^*_{\R^n} \sup \left\{ \prod_{i=1}^m f_i^{c_i}(y_i)\,; \;
	\sum_{i=1}^m c_i B_i^* y_i=x \mbox{ and } y_i\in \R^{n_i} \right\}\, dx ,$$
and
$$ J((f_i)_{i=1}^m)= \int_{\R^n}  \prod_{i=1}^m f_i^{c_i}(B_i x) \, dx.$$
Let $E$ be the largest constant such that for all $(f_i)_{i=1}^m$,
$${}\hspace{4 cm} I((f_i)_{i=1}^m) \ge E \prod_{i=1}^m \left( \int_{\R^{n_i}} f_i \right) ^{c_i},
	 \hspace{4 cm} (\mathrm{RBL})$$
and let $F$ be the smallest one such that for all $(f_i)_{i=1}^m$,
$${}\hspace{4 cm} J((f_i)_{i=1}^m) \le F \prod_{i=1}^m \left( \int_{\R^{n_i}} f_i \right) ^{c_i}.
	 \hspace{4 cm} (\mathrm{BL})$$
Then $E$ and $F$ can be computed just with centered Gaussian functions, that is
$$E =\inf \left\{ \frac{I((g_i)_{i=1}^m)}{ \prod_{i=1}^m \left( \int_{\R^{n_i}} g_i \right) ^{c_i}} \,; \;
		g_i \mbox{ centered Gaussian on } \R^{n_i}, \, i=1, \ldots,m \right\},$$
and
$$F =\sup \left\{ \frac{J((g_i)_{i=1}^m)}{ \prod_{i=1}^m \left( \int_{\R^{n_i}} g_i \right) ^{c_i}} \,; \;
		g_i \mbox{ centered Gaussian on } \R^{n_i}, \, i=1, \ldots,m \right\},$$
Moreover, if we denote by $D$ the largest real number such that
$$\det \left(\sum _{i=1}^m c_i B^*_i A_i B_i \right) \ge D \prod _{i=1}^m \left( \det A_i \right)^{c_i},$$
for all $A_i \in \mathcal S^+ (\R^{n_i}), \, i=1,\ldots,m$, then
$$ E=\sqrt D \qquad \qquad \mbox{ and } \qquad \qquad F=\frac1{\sqrt D} \cdot $$

\end{theo}
{\bf Remark 1:} the hypothesis $\sum_{i=1}^m c_i n_i =n$ is just a necessary homogeneity
	condition for $E$ to be positive and for $F$ to be finite. The condition on $\cap \ker B_i$
	ensures that $\sum_{i=1}^m c_i B^*_i A_i B_i$ is an isomorphism. Actually, the conclusion
	of the theorem remains valid without this condition, but it is obvious because $D=0$. 

\noindent
{\bf Remark 2:} Notice that the reverse Brascamp-Lieb inequality for $m=2$, $n_1=n_2=n$,
	$B_1=B_2=B_1^*=B_2^*=I_n$ and $c_1=\alpha=1-c_2$, where $I_n$ is the indentity map
	on $\R^n$ and $0<\alpha <1$, is the inequality of Pr\'ekopa-Leindler. Indeed the
	constant $D$ is
	$$ D= \inf\limits_{A_1,A_2 \in \mathcal S ^+(\R^n)} \frac{\det (\alpha A_1 + (1-\alpha)A_2)}
		{(\det A_1)^{\alpha} (\det A_2)^{1-\alpha}}=1$$
	by the arithmetic-geometric inequality. So (RBL) becomes, for all $f,g \in L_1^+(\R^n)$,
	$$  \int\limits_{\R^n}^* \sup\limits_{x=\alpha u + (1-\alpha)v } 
		f^{\alpha}(u) g^{1-\alpha}(v) \, d^n x 
		\ge \left(\int_{\R^n} f \right)^{\alpha} \left(\int_{\R^n} g \right)^{1-\alpha}. $$
	It is well-known that this inequality implies the celebrated Brunn-Minkowski theorem: for
	$A,B$ compact non-void subsets of $\R^n$,
	$$ \mathrm{Vol}^{\frac1n} (A+B) \ge  \mathrm{Vol}^{\frac1n} (A) + \mathrm{Vol}^{\frac1n} (B)\, .$$

\vspace{5 mm}
The proof of Theorem~\ref{thRBL} is divided into lemmas.
	We deal first with the study of the behaviour of $I$ and $J$
 	with respect to centered Gaussian functions. We set
$$E_g =\inf \left\{ \frac{I((g_i)_{i=1}^m)}{ \prod_{i=1}^m \left( \int_{\R^{n_i}} g_i \right) ^{c_i}}; \;
		g_i \mbox{ centered Gaussian on } \R^{n_i}, \, i=1, \ldots,m \right\},$$
and
$$F_g =\sup \left\{ \frac{J((g_i)_{i=1}^m)}{ \prod_{i=1}^m \left( \int_{\R^{n_i}} g_i \right) ^{c_i}}; \;
		g_i \mbox{ centered Gaussian on } \R^{n_i}, \, i=1, \ldots,m \right\},$$
our aim is to prove that $E=E_g=\sqrt D$ and $F=F_g=D^{-1/2}.$
We begin by a classical computation, done in \cite{brasl76bcyi}; its only uses the fact that if
$M\in \mathcal S^+ (\R^k)$ then
$$ \int_{\R^k} \exp(-\langle x,Mx \rangle)\, dx = \sqrt{ \frac{\pi^k}{\det M}} \cdot $$

\begin{lem}
\label{gaussbl}
	With the notations of Theorem~\ref{thRBL}, we have
	$$ F_g=\frac1{\sqrt D} \cdot$$
\end{lem}
\option{
{\bf Proof:} For $i=1, \ldots,m$, let $A_i \in \mathcal S^+ (\R^{n_i})$ and let $g_i$ be defined
	on $\R^{n_i}$ by $$g_i(x)=\exp(-\langle A_ix,x \rangle),$$ then
	\begin{eqnarray*}
		J(g_1,\ldots,g_m)&=& \int_{\R^n} \exp \left( - \langle
			 \sum_{i=1}^m c_i B_i^* A_i B_i x, x \rangle \right) \, dx \\ 
		 &=& \left( \frac{\pi^n}{\det \left(\sum_{i=1}^m c_i B_i^* A_i B_i \right)} \right)^{1/2}.
	\end{eqnarray*}
	As
	$$ \prod_{i=1}^m \left( \int_{\R^{n_i}} g_i \right) ^{c_i}
		= \prod_{i=1}^m \left( \frac{\pi^{n_i}}{\det A_i} \right)^{c_i/2},$$
	we obtain
	$$ \frac{J(g_1, \ldots ,g_m)}{ \prod_{i=1}^m \left( \int_{\R^{n_i}} g_i \right) ^{c_i}}
		=\left( \frac{\det \left(\sum_{i=1}^m c_i B_i^* A_i B_i \right)}
		     {\	\prod_{i=1}^m \left( \det A_i \right)^{c_i}}
		 \right) ^{-1/2}.$$
	The result of the lemma follows by taking the supremum. \qed
	}

Our next lemma links $E_g$ and $F_g$ by means of duality between quadratic forms.
\begin{lem}
\label{gaussrbl}
	 With the previous notations, we have
	$$E_g \cdot F_g=1\, ,$$
	and $E_g=0$ if and only if $F_g=+\infty$.
\end{lem}
{\bf Proof:} For $i=1, \ldots,m$, let $A_i \in \mathcal S^+ (\R^{n_i})$ and let $Q$
	be the quadratic form on $\R^n$ defined by
	$$Q(y)=\langle \sum_{i=1}^m c_i B_i^* A_i B_i y, y \rangle.$$
	Let $Q^*$ be the dual quadratic form of $Q$, we recall that it is defined on $\R^n$ by
	$$Q^*(x)= \sup \left\{ |\langle x,y \rangle|^2 \, ; \; Q(y) \le 1\right\}.$$
	We also intoduce the application $R$ on $\R^n$ such that for all $x \in \R^n$,
	$$R(x)=\inf \left\{ \sum_{i=1}^m c_i \langle A_i^{-1} x_i,x_i \rangle \, ; \;
		x=\sum_{i=1}^m c_i B^*_i x_i \mbox{ and for all $i$, } x_i\in \R^{n_i} \right\}.$$
	We show now that $R=Q^*$.
	Indeed, assume that $x=\sum_{i=1}^m c_i B^*_i x_i$ with $x_i\in \R^{n_1}$ for $i=1, \ldots, m$,
	then
	$$|\langle x,y \rangle|^2 =  | \langle \sum_{i=1}^m c_i B^*_i x_i,y \rangle |^2 
		= | \sum_{i=1}^m \langle \sqrt{c_i} A_i^{-1/2} x_i , \sqrt{c_i} A_i^{1/2} B_i y \rangle|^2.
	$$
	By the Cauchy-Schwartz inequality, applied to the quadratic form $\phi$ on $\R^{n_1}\times \cdots
	\times \R^{n_m}$ defined by $\phi(X_i,\ldots,X_m)= \sum_{i=1}^m \langle X_i,X_i \rangle$, one gets:
	\beq	|\langle x,y \rangle|^2
			&\le& \left(  \sum_{i=1}^m |\sqrt{c_i} A_i^{-1/2} x_i |^2 \right)
				\left( \sum_{i=1}^m |\sqrt{c_i} A_i^{1/2} B_i y |^2 \right) \\
			&=& \left(  \sum_{i=1}^m c_i \langle x_i,A_i^{-1} x_i \rangle \right)
				\left( \langle \sum_{i=1}^m c_i B^*_i A_i B_i y , y \rangle \right).
	\eeq
	In fact, one easily checks that there is equality in the previous argument if one takes
	$$ y= \left( \sum_{i=1}^m c_i B^*_i A_i B_i \right)^{\!\! -1} \! \! x$$
	and $$x_i =A_i B_i y \qquad i=1, \ldots, m ,$$
	therefore $R=Q^*$.

	We apply this result to our integrals of Gaussian functions. Straightforward computations
	give that
	$$\frac{J(G_{A_1},\ldots, G_{A_m})}{\prod _{i=1}^m \left(\int G_{A_i}\right) ^{c_i}}
		= \sqrt{ \frac{ \prod _{i=1}^m ( \det A_i)^{c_i}}{\det Q}}, $$
	and
	$$\frac{I(G_{A_1^{-1}},\ldots, G_{A_m^{-1}})}{\prod _{i=1}^m \left(\int G_{A_i^{-1}}\right) ^{c_i}}
		= \sqrt{ \frac{ \prod _{i=1}^m ( \det A_i)^{-c_i}}{\det R}} \cdot $$
	Using the result $R=Q^*$ and the classical duality relation $\det Q \cdot \det Q^* =1$, one has
	$$\frac{J(G_{A_1},\ldots, G_{A_m})}{\prod _{i=1}^m \left(\int G_{A_i}\right) ^{c_i}}
		\cdot
	  \frac{I(G_{A_1^{-1}},\ldots, G_{A_m^{-1}})}{\prod _{i=1}^m \left(\int G_{A_i^{-1}}\right) ^{c_i}}=1,$$
	therefore $E_g = F_g^{-1}$. \qed
\vspace{4 mm}

\noindent
{\bf Remark:} We emphasize the equivalence for $A_i \in \mathcal S ^+(\R^n)$, $i=1,
	\ldots, m$, of the assertions
	\begin{itemize}
		\item $\D \det \left(\sum _{i=1}^m c_i B^*_i A_i B_i \right) 
			=D \prod _{i=1}^m \left( \det A_i \right)^{c_i}.$
		\item The $m$-tuple of centered Gaussians $(G_{A_1}, \ldots, G_{A_m})$
			is a maximizer for (BL).
		\item The $m$-tuple of centered Gaussians $(G_{A_1^{-1}}, \ldots, G_{A_m^{-1}})$
			is a minimizer for (RBL).
	\end{itemize}

\vspace{4 mm}
We state now the fundammental result which, combined with the two previous lemma, will
	suffice to establish Theorem~\ref{thRBL}. Since the theorem is already established if
	$D=0$, we assume from now on that $D$ is positive.
\begin{lem} 
\label{fond}
	For $i=1...m$, let $f_i$ and $h_i$ belong to $L^+_1(\R^{n_i})$ and satisfy 
	$\int_{\R^{n_i}} f_i = \int_{\R^{n_i}} h_i.$ Then 
	$$ I(f_1, \ldots, f_m) \ge D \cdot J(h_1, \ldots, h_m). $$
\end{lem}
In \cite{bart97iblc}, the author proved this result for functions of one real variable,
using measure-preserving mappings; given $f$ and $h$, two non-negative functions on $\R$
with integral one, there exists a non-decreasing mapping $u$ such that for all $x \in \R$:
	$$ \int_{- \infty}^{u(x)} f = \int_{- \infty}^x h \, .$$
In other words, $u$ maps the probability measure of density $h$ onto the probability
measure of density $f$. Our proof in the general case (i.e. for functions of several
variables) is also based on measure-preserving mappings. But, in dimension larger than one,
there is a large choice of such mappings between two sufficiently regular
probability measures. For our purpose, the Brenier mapping (see \cite{bren87dprm},
 \cite{bren91pfmr}) fits perfectly; it has the additionnal convenient property
of deriving from a convex potential. Brenier proved its existence and uniqueness
under certain integrability assumptions on the moments of the measures, which where
later removed by McCann \cite{mcca94ctig}, \cite{mcca95eumm}. Let us state the result
that we need.
\begin{theo}
	Let $f_1,f_2$ be non-negative measurable functions on $\R^n$ with integral one. There 
 	exists a convex function $\phi$ on $\R^n$ such that the map $u= \nabla \phi$ has 
	the following property: for every non-negative borelian function $b$ on $\R^n$,
	$$ \int_{\R^n} b(u(x)) f_2(x) \, dx = \int_{\R^n} b(x) f_1(x) \, dx. $$
\end{theo}
The function $\phi$ given by this theorem represents a generalized solution
of the Monge-Amp\`ere equation
	$$ \det ( \nabla^2 \phi(x)) f_2(\nabla \phi(x)) = f_1(x). $$
In fact, the gradient of $\phi$ is unique $f_1 \, dx$-almost everywhere.
Since it will be convenient to work with strong solutions, we recall here a corollary
of a theorem of Caffarelli \cite{caff92rmcp}, who has developped a regularity theory for
these convex solutions. 
\begin{theo}
	For $i=1,2$, let $\Omega_i$ be bounded domains of $\R^n$ and let $f_i$ be non-negative
	functions, supported on  $\Omega_i$. Assume that $f_i$ and $1/f_i$ are bounded on 
	$\Omega_i$ and that $\Omega_2$ is convex. If $f_i$, $i=1,2$ are Lipschitz then the
	Brenier mapping $\phi$ is twice continuously differentiable.
\end{theo}
Let $\mathcal C _L(\R^n)$ be the set of functions $f \in L_1^+(\R^n)$ which are the restriction
	to some opened Euclidean ball of a positive Lipschitz function on $\R^n$.

Let us remark that it suffices to establish (BL) and (RBL) for functions in  $\mathcal C _L(\R^{n_i})$.
	We use strongly the monotonicity of the applications $I$ and $J$.
	By the regularity of measure, for any $f \in  L_1^+(\R^n)$ and any
	$\varepsilon > 0$, there exists a function $s$, which is a positive combination
	of characteristic functions of compact sets, such that 
	$$ f \ge s \mbox{  and  } \int f - \int s \le \varepsilon \, ,$$
	so its is enough to prove (RBL) for such functions. As $s$ is clearly the
	pointwise limit of some decreasing sequence of Lipschitz functions, it suffices
	to work on Lipschitz functions. Moreover, we can assume
	these functions to be positive (by adding some Gaussian $G/N$, where $N$ tends to
	infinity).
	Eventually, by truncation, it is enough to work with functions
	in $\mathcal C _L(\R^{n})$.

The same kind of argument is valid for (BL). Moreover,
	since (BL) is equivalent to the boundedness of a multilinear operator which is, with
	respect to each function, a linear kernel operator (because the $B_i$'s are surjective)
	, it clearly suffices to show (BL) for a dense subset of $L_1$.

{\bf Proof of Lemma~\ref{fond}: }
	By homogeneity we can assume that $\int f_i =\int h_i =1$ for all $i$. The previous
	remark allows us to work with functions $f_i,h_i$ belonging to $\mathcal C _L(\R^{n_i})$,
	so that we can use Caffarelli's regularity result and Brenier
	theorem. We denote by $\Omega_{h_i}$ the domain where $h_i$ is positive. 
	We get, for $i=1,\ldots, m$, differentiable mappings $T_i$ deriving
	from convex potentials and such that, for all $x \in \Omega_{h_i}$,
	$$ \det \left( dT_i(x) \right) \cdot f_i(T_i(x)) = h_i (x).$$
	Since $T_i$ derives from convex potential, its differential is symmetric semi-definite positive
	and because of the previous equation and of the non-vanishing property of $h_i$,
	we know that for all $x \in \Omega_{h_i}$, $dT_i(x) \in \mathcal S^+ (\R^{n_i})$.
	
	We define a function $\Theta$ from $\cap_{i=1}^m B_i^{-1}(\Omega_{ h_i}) \subset \R^n$
	into $\R^n$ by
	$$ \Theta (y) = \sum_{i=1}^m c_i B_i^* ( T_i ( B_i y)) \,.$$
	Its differential is symmetric semi-definite positive
	$$ d \Theta (y) = \sum_{i=1}^m c_i B_i^* dT_i (B_i y) B_i \,,$$
	and it is actually definite positive because:
	$$ \det \left(\sum_{i=1}^m c_i B_i^* dT_i (B_i y) B_i \right) \ge D
		\prod_{i=1}^m  \left( \det dT_i(B_i y) \right) ^{c_i} > 0 $$
	In particular for all $v \neq 0$ of $\R^n$,
	$$ \langle d\Theta (y) \cdot v , v\rangle  > 0 $$
	so $\Theta$ is injective.
	Denoting $S=\cap_{i=1}^m B_i^{-1}(\Omega_{ h_i})$, we can write
	\beq
		\int\limits_{\R^n} \prod_{i=1}^m h_i^{c_i} (B_i y) \, dy &=& 
			\int\limits_{S} \prod_{i=1}^m h_i^{c_i} (B_i y) \, dy\\
		&=& \int\limits_{S} \prod_{i=1}^m \left(
			f_i(T_i(B_i y)) \det dT_i(B_i y) \right)^{c_i} \, dy \\
		& \le& \frac1D  \int\limits_{S} 
			\prod_{i=1}^m f_i(T_i(B_i y))^{c_i}
			\det \left(\sum_{i=1}^m c_i B_i^* dT_i (B_i y) B_i \right) \, dy \\
		& \le& \frac1D \int\limits_{S}
			 \sup\limits_{\Theta(y)=\sum_{i=1}^m c_i B_i^* x_i}
			\left(\prod_{i=1}^m f_i(x_i)^{c_i} \!\right)
			\det(d\Theta(y)) \, dy \\
		& \le& \frac1D \int\limits_{\R^n} \sup\limits_{x=\sum_{i=1}^m c_i B_i^* x_i}
			\left( \prod_{i=1}^m f_i(x_i)^{c_i} \! \right)
			\, dx \,
	\eeq
	which concludes the proof.
	\qed

\section{Equality cases}
In this section, we restrict to functions of one real variable.
	With the notations of Theorem~\ref{thRBL} there are vectors $v_1, \ldots, v_m$, in $\R^n$
	such that span$((v_i)_{i=1}^m)= \R^n$ and for all $x \in \R^n$ and $t \in \R$,
	$$ B_i(x)= \lag x, v_i \rag \mbox{   and   } B_i^*(t) = t v_i \, .$$
	We are going to study the best constant in inequalities (BL) and (RBL)
	and to characterize equality cases. We call maximizers the non-zero functions
	that give equality in (BL) and minimizers those that provide equality
	in (RBL).

\subsection{The geometric structure of the problem}

We introduce some notations. For a subset $K$ of $\{1,\ldots, m\}$, we denote by
$E_K$ the linear span in $\R^n$ of the vectors $(v_k)_{k \in K}$. We will call
{\it adapted partition} a partition $S$ of $\{1,\ldots, m\}$ such that:
$$ \R^n = \bigoplus_{K \in S} E_K \; .$$
These partitions are usefull because\option{, as we will see later,} this splitting
of the space $\R^n$ implies a splitting of the Brascamp-Lieb inequality
 and of its converse, so that one just needs to work separately on each piece.
We shall show first that there exists a best adapted partition.
\begin{prop}
\label{decomp}
	Let $\bowtie$ be the relation on $N_m=\{1,\ldots, m\}$ defined by as follows:
	$ i \bowtie j$ if and only if there exists a subset $K$ of $N_m$ of 
	cardinality $n-1$ such that both $(v_i, (v_k)_{k \in K})$ and
	 $(v_j, (v_k)_{k \in K}) $ are basis of $\R^n$. Let $\sim$ be the
	transitive completion of $\bowtie$ ($i \sim j$ means that there exits
	a path between $i$ and $j$ in which two consecutive elements are in relation
	for $\bowtie$).

	Then $\sim$ is an equivalence relation and the subdivision $C$
	of $N_m$ into equivalence classes for $\sim$ is the most accurate adapted
	partition.
\end{prop}
{\bf Proof:} We establish first that $C$ is more accurate than any adapted partition $S$.
	Let $I,J \in S$, $I \neq J$ and let $i \in I$, $j \in J$. It suffices to show that 
	$ i \bowtie j$ is impossible.

	Assume precisely that $ i \bowtie j$, there exists $K \subset N_m$ such that
	$$ \mathbf e _i=(v_i, (v_k)_{k\in K}) \mbox{ and } \mathbf e _j=(v_j, (v_k)_{k\in K})$$
	form basis of $\R^n$. As $S$ is adapted, we have $\R^n = \bigoplus_{H \in S} E_H$,
	each of them being spanned by certain $v_i$'s. So, every basis of $\R^n$ which
	is formed of some of the $v_i$'s must contain $\dim (E_H)$ elements in $E_H$.
	But our basis $  \mathbf e _i$ and $  \mathbf e _i$ do not have the same number
	of vectors in $E_I$ because $v_i \in E_I$ and $v_j \in E_J$. Thus we have a contradiction.

	\vspace{4 mm}

	We prove now that the partition $C$ is adapted to our geometric setting. Let $I$ be an
	equivalence class for $\sim$ and let $E_I$ be the corresponding space. Since 
	the vectors $(v_i)_{i=1}^m$ span all $\R^n$, we find a permutation of indices
	such that $\mathbf b =(v_1, \ldots, v_n)$ is a basis of $\R^n$ and $(v_1, \ldots, v_r)$
	is a basis of $E_I$ for a certain $r \le n$. Let us denote by $F$ the span
	of $v_{r+1}, \ldots , v_n$.

	Let $i \in N_m$; the vector $v_i$ can be decomposed in the basis $\mathbf b$:
	$$ v_i =\sum_{i=1}^n \alpha_i v_i \, .$$
	For any $j \le n$, we notice
	$$ \det_{\mathbf b} ( v_1, \ldots, v_{j-1}, v_i, v_{j+1}, \ldots , v_n )= \alpha_j \, ,$$
	hence $\alpha_j \neq 0$ implies that $v_i$ and $v_j$ belong to neighbourg
	basis, that is $i \sim j$.
	So, if $i \in I$, as $i$ can be in relation for $\bowtie$ only with elements
	of $I$ we have $\alpha_{r+1}, \ldots, \alpha_n=0$. Thus $i \in I$ implies 
	$v_i \in E_I$.
	By a similar argument, if $i \not \in I$, $\alpha_1, \ldots, \alpha_r=0$
	and $v_i$ belongs to $F$.
	We have proved that $\R^n= \mbox{span}\{ v_i, i \in I\} \bigoplus 
		\mbox{span}\{ v_i, i \not \in I\}$, this is the first step
	of the decomposition. The result follows by induction, noticing that
	the relation $\sim$ can be restricted to $F$. \qed

\vspace{5 mm}
\option{
{\bf Principle:} {\it  Adapted decompositions are convenient because the Brascamp-Lieb inequality
	on $\R^n$ is equivalent to the corresponding inequalities on 
	the subspaces which appear in the decomposition, moreover there
	is equality in the inequality (for non identically zero functions)
	if and only if there is equality on each subspace.}

The same principle is valid for the reverse inequality (but we will not prove it because
	the proof is similar and straightforward).

{\bf Proof:} Assume that $S$ is an adapted decomposition. With the previous notations
	we have $\R^n= \bigoplus_{I \in S} E_I$. Let $n_I$ be the dimension
	of $E_I$. For each $I$ let $u_I$ be a linear isomorphism between $E_I$
	and $\R^{n_I}$ and let $u=\oplus_{I \in S} u_I$ be an isomorphism of
	$\R^n$  such that the images of two different subspaces $E_I, E_J$ are
	orthogonal in $\R^n \simeq \oplus_{I \in S}^{\bot} \R^{n_I}$. For 
	$y \in \R^n$ we will denote by $(y_I)_{I \in S}$ its coordinates
	in the previous orthogonal sum.
	This map rectifies the direct sum and will allow us to apply Fubini's
	theorem. We can now illustrate the splitting of the quantities involved
	in the Brascamp-Lieb inequality:
	\beq
		 \int_{\R^n} \prod_{i=1}^m f_i^{c_i} ( \langle x, v_i \rangle) \, d^n x
			&=& \det(u) \int_{\R^n} \prod_{i=1}^m f_i^{c_i} ( \langle y, uv_i \rangle) \, d^n x
			 \\
		&=&   \det(u) \int_{\R^n} \prod_{I \in S} \left( \prod_{i\in I} 
			f_i^{c_i} ( \langle y, uv_i \rangle)  \right)\, \prod_{I \in S} d^{n_I} y_I
	\eeq

	Using the orthogonality properties of $u$ we have, when $i \in I$,
	$$\langle y, uv_i \rangle= \sum_{J \in S} \langle y_J, uv_i \rangle
			=\langle y_I, uv_i \rangle \, .$$
	Thus
	\beq 
		 \int_{\R^n} \prod_{i=1}^m f_i^{c_i} ( \langle x, v_i \rangle) \, d^n x  
		&=& \prod_{I \in S}  \left( \int_{\R_{n_I}} \prod_{i\in I} f_i^{c_i} (\langle y_I, uv_i \rangle)
			d^{n_I}y_I \cdot \det(u_I) \right) \\
		&=& \prod_{I \in S} \left(\int_{E_I} \prod_{i\in I} f_i^{c_i} (\langle y, v_i \rangle)
			d^{E_I} y \right)\, .
	\eeq
	On the other hand its is obvious from the partition property that
	$$ \prod_{i=1}^m \left( \int f_i \right)^{c_i} =\prod_{I \in S} \left( 
		\prod_{i\in I} \left( \int f_i \right)^{c_i} \right) \, .$$
	\qed
\vspace{4 mm}
}

\noindent
As a conclusion let us notice that is suffices to study the case when the relation $\bowtie$
	has only one equivalence class. In this case we say that $(\R^n, (v_i)_{i=1}^m)$
	is {\it irreducible}.

\subsection{The Gaussian case}
Let $v_1, \ldots, v_m$ be the vectors of $\R^n$ defined at the beginning of this section.
 For $I \subset \{1, \ldots, m \}$ of 
cardinal $|I|=n$, we denote
	$$ d_I =  \det( (v_i)_{i \in I})^2 .$$
For each $m$-tuple $c=(c_i)_{i=1}^m$ of positive real, we study the constant $D_c$ defined by
	$$ D_c= \inf \left\{ \frac{ \det (\sum_{i=1}^m \lambda_i v_i \otimes v_i)}
				{ \prod_{i=1}^m \lambda_i^{c_i}}   \, ; \;
			  \lambda_i >0, \, i=1 \ldots m
		\right\}. $$
We wish to know when it is positive and when it is achieved. We will sometimes call minimizers 
	the $m$-tuples $(\la_i)_{i=1}^m$ for which $D_c$ is achieved.

The computation of the previous determinant is made possible by the Cauchy-Binet
	formula which we recall:
\begin{prop}
	Let $m \ge n$ be integers; let $A$ be a $n \times m$ matrix and let $B$ 
	be a $m \times n$ matrix. For $I \subset N_m$ of cardinality $n$ we denote
	by $A_I$ the square matrix obtained from $A$ by keeping only the columns 
	with indices in $I$; we denote by $B^I$ the square matrix obtained from
	$B$ by keeping the rows with indices in $I$. Then we have the formula
	$$ \det (AB) = \sum_{|I|=n} \det (A_I) \det (B^I)$$
	where the sum is over the subsets of cardinality $n$ of $N_m$.
\end{prop}
The relevance of this formula is clear from 
\begin{coro}
\label{cb}
	Let $m \ge n$ and let $(v_1, \ldots , v_m)$ be vectors in $\R^n$, then
	$$ \det\left( \sum_{i=1}^m \lambda_i v_i \otimes v_i \right)
		= \sum_{|I|=n} \lambda_I \left( \det((v_i)_{i \in I}) \right)^2, $$
	where for $I \subset N_m$, we have set $\lambda_I = \prod_{i \in I} \lambda_i .$
\end{coro}

The condition for $D_c$ to be non-zero has a rather nice geometric expression which
requires some notations. For $ I \subset \{1, \ldots, m \}$, we denote by
$1_I$ the vector of $\R^m$ of coordinates $(1_I)_i= \delta_{i \in I}$ (it is
the characteristic function of $I$). We denote by $c$ the vector $(c_i)_{i=1}^m$.
One has the following result:
\begin{prop}
\label{Dpas0}
	The infimum $D_c$ is positive if and only if the vector $c$ belongs to the
	convex hull of the characteristic vectors $1_I$ of the subsets $I$ of
        cardinal $n$ such that the vectors $ (v_i)_{i \in I}$ form a basis of $\R^n$.
\end{prop}
{\bf Proof:}
	We shall show first that the condition is sufficient. Assume that we have
	a family of non-negative real numbers $(t_I)_{|I|=n}$ indexed by the subsets
	of cardinal $n$ of $\{1, \ldots, m\}$, such that
	$$ \begin{array}{l} t_I=0 \mbox{ whenever } d_I=0, \\
			c_i=\sum_{|I|=n, \, i\in I} t_I, \, \mbox{ for all } i
	   \end{array} $$
	Let $\lambda_i$, $i=1, \ldots, m$ be positive. By the Cauchy-Binet formula,
	we have:
	$$\det (\sum_{i=1}^m \lambda_i v_i \otimes v_i)= \sum_{|I|=n} \lambda_I d_I
		=\sum_{t_I \neq 0} t_I \left( \frac{\lambda_I}{t_I} d_I \right) 
			+  \sum_{t_I=0} \lambda_I d_I.
	$$
	The second term is non-negative. We apply the arithmetic mean-geometric mean inequality
	with coefficients $t_I$ (their sum is indeed one), and for each $i$ we gather the
	factors with $\lambda_i$. Each $\lambda_i$ will appear with an exponent equal to
	$$ \sum_{i \in I, \, t_I \neq 0 } t_I,$$
	this is $c_i$ by hypothesis. Thus we have
	$$\det (\sum_{i=1}^m \lambda_i v_i \otimes v_i) \ge 
		\prod_{t_I \neq 0} \left(\frac{d_I}{t_I}\right)^{t_I} \prod_{i=1}^m \lambda_i^{c_i}.$$
	Since $t_I \neq 0 $ implies $d_I \neq 0$, the constant $D_c$ is positive. 

	Let us prove now that the condition is necessary. Assume that the function
	$$\Delta (\lambda_1, \ldots, \lambda_m) =\frac { \det (\sum_{i=1}^m \lambda_i v_i \otimes v_i)}
							{ \prod_{i=1}^m \lambda_i^{c_i}} $$
	is bounded below by a positive $D_c$ when $\lambda_i >0, \, i=1 \ldots m$.
	If we take $\lambda_i = N^{-x_i}$, where the $(x_i)$ are arbitrary, and 
	$N$ tends to the infinity, then $\Delta (\lambda_1, \ldots, \lambda_m)$ is
	equivalent to a positive constant times $N$ to the exponent:
	$$\sum_{i=1}^m x_i c_i +\max \{ -\sum_{i \in I} x_i \, ; \; |I|=n \mbox{ and } d_I \neq 0 \}.$$
	As by hypothesis, $\Delta$ cannot tend to 0, the exponent must be non-negative. Thus for 
	all $(x_i)_{i=1}^m \in \R^m$, one has:
	$$ \sum_{i=1}^m x_i c_i \ge
		 \min \left\{ \sum_{i \in I} x_i \, ; \; |I|=n \mbox{ and } d_I \neq 0 \right\}.$$
	Equivalently, for all $x \in \R^m$,
	$$ \min_{ |I|=n, \, d_I \neq 0} \langle x, 1_I \rangle  \le  \langle x, c \rangle,$$
	which can be reformulated in terms of convex cones as	
	$ \cap_{d_I \neq 0} \mathcal C _{1_I} \subset \mathcal C _c \,$,
	where, for $y \in \R^m$, $\mathcal C _y = \{ x \in \R^m \, ; \;  \langle x, y \rangle \ge 0 \}.$
	By duality of convex cones, this implies that the vector $c$ belongs to the
	convex cone generated by the vectors $1_I$ such that $d_I \neq 0$.
	Thus there exist non-negative real numbers $(t_I)_{I, d_I \neq 0}$ such
	that for all $i\le m$,
	$$c_i=\sum_{|I|=n \mbox{ and } i\in I} t_I.$$
	If we make the sum on $i$ of the previous relations, we get $\sum_{d_I \neq 0} t_I=
	(\sum_{i=1}^m c_i)/n$. But the Hypothesis $D_c>0$ implies that the numerator and the
	denominator of $\Delta$ must be of the same homogeneity degree in the variables, so
	$\sum_{i=1}^m c_i=n$ and we have derived that $c$ belongs to the convex hull of 
	the $1_I$ such that $d_I \neq 0$. \qed

{\bf Remark:} Let $K= \{ x\in [0,1]^m \, ;\;  \sum_{i=1}^m x_i=n \}$, it is the convex hull
	of the vectors $(1_I)_{|I|=n}$. By the previous result, $D_c$ is non-zero only if
	$c$ is in $K$. If the vectors $(v_i)$ are in generic position, $D_c \neq 0$ if and only
	if $c \in K$. But as the $1_I$ are clearly the only extremal points of $K$, any geometrical
	degeneracy (i.e. any $d_I$ equal to zero) will imply a reduction of the domain where
	$c$ must be.

We know that $D_c$ is positive if and only if $c$ can be written as a convex combination
	of certain vectors. The next proposition states that $D_c$ is achieved if and only
	if there exists a convex combination with some additional property.
\begin{prop}
	The constant $D_c$ is achieved if and only if there exist positive numbers $(t_I)_{|I|=n}$
	and $(\lambda_i)_{i=1}^m$ such that
	$$ c= \sum_{|I|=n} t_I 1_I$$
	and for all $I$ 
	$$ t_I = d_I \prod_{i \in I} \lambda_i \, .$$
\end{prop}
Notice that $d_I=0$ implies $t_I=0$, so the result is coherent with the previous one.

{\bf Proof:} 
	The {\em if} part comes from a precise study of the proof of proposition~\ref{Dpas0}:
	the inequality it gives is an equality for the $m$-tuple $(\lambda_i )$ because
	for all $I$, $\lambda_I d_I / t_I=1$ so the arithmetic-geometric inequality
	is an equality; moreover the term $\sum_{I, \, t_I=0} \lambda_I d_I$ is zero.
	\option[The {\em only if} part is obvious by differentiation.]
	{ 
	Assume that $D_c$ is achieved for a $m$-tuple $(\mu_i)_{i=1}^m$; the function $\Delta$ 
	has a minimum at this point, so its differential must be zero. A simple computation
	implies that for all $i \le m$, one has
	$$ c_i = \frac{ \sum_{I, \, i \in I} \mu_I d_I}{\sum_I \mu_I d_I} \cdot $$
	So,
	$$ c= \sum_{I} \left( \frac{ \prod_{i \in I} \mu_i}{\sum_I \mu_I d_I} d_I \right)
		1_I \, ,$$
	and the coefficients have the property that we predicted.
	}
	\qed

\vspace{3 mm}
We are going to rewrite our problem in the setting of Fenchel duality for convex functions
in order to use the following result (see \cite{rockCA} p264):
\begin{prop}
\label{dualfenchel}
	Let $\phi$ be a l.s.c. convex function on $\R^m$ and let $\phi^*$ be its Fenchel conjugate,
	defined for $x \in \R^m$ by
	$$ \phi^* (x)= \sup_{y \in \R^m} \langle x,y \rangle - \phi(y) \, .$$
	Then $\phi^*(x)$, which is a supremum, is achieved if and only if
	$\phi^*$ is subdifferentiable at the point $x$. In particular, it is 
	achieved when $x$ belongs to the relative interior of $\mathrm{dom}(\phi^*)=
	\{ y \in \R^m \, ; \; \phi^*(y) < + \infty \} $.
\end{prop}
Let us define the function $\phi$ on $\R^m$ by
	$$ \phi( x_1, \ldots, x_m)= \log \det \left( \sum_{i=1}^m e^{t_i} v_i \otimes v_i \right),$$
	the next proposition links our problem on $D$ with the study of the Fenchel conjugate
	of $\phi$.

\begin{prop}
\label{lienDphi}
	\begin{enumerate}
		\item The function $\phi$ is convex.
		\item The constant $D_c$ is equal to $\exp( - \phi^* (c))$.
		\item $D_c$ is positive if and only if $c \in \mathrm{dom}(\phi^*)$.
		\item $D_c$ is achieved if and only if $\phi^*(c) $ is.
		\item $\mathrm{dom}(\phi^*)$ is equal to $K=\mathrm{conv}\{1_I \, ; \; d_I \neq 0 \}$.
		\item The constant $D_c$ is achieved when $c$ belongs to the relative interior
			of $K$.
	\end{enumerate}
\end{prop}
{\bf Proof:}
	The convexity of $\phi$ is a consequence of the Cauchy-Schwartz inequality: let $s,t \in \R^m$,
	\beq 
		\phi \left(\frac{t+s}{2} \right) \option{&=& \log \left( \sum_{|I|=n} \exp \left(\sum_{i \in I}
			(t_i+s_i)/2 \right) d_I \right) \\}
		&=&\log \left( \sum_{|I|=n} \left\{d_I \exp \left(\sum_{i \in I} t_i \right)\right\}^{\frac12}
			\left\{d_I \exp\left(\sum_{i \in I} s_i \right)\right\}^{\frac12} \right) \\
		&\le& \log \left( \left\{ \sum_{i \in I}d_I \exp\left(\sum_{i \in I} t_i\right)\right\}^{\frac12}
				\left\{\sum_{i \in I}d_I \exp\left(\sum_{i \in I} 
						s_i\right)\right\}^{\frac12}  \right) \\
		&=& \frac{\phi(t)+\phi(s)}{2} \cdot
	\eeq
	\option[The other assertions are also very simple.]
	{
	The link between $D_c$ and $\phi^*$ is straightforward:
	\beq
		D_c &=& \inf_{\lambda_i >0} \frac{\det \left(\sum_{i=1}^m \lambda_i v_i \otimes v_i\right)}
			{\prod_{i=1}^m \lambda_i^{c_i}} \\
		&=& \exp \left( - \sup_{x_i \in \R} \left\{
			-\log \left( \sum_{i=1}^m e^{x_i} v_i \otimes v_i \right)
			+ \sum_{i=1}^m c_i x_i  \right\} \right) \\
		&=& \exp \left( - \phi^*(c) \right).
	\eeq
	The other conclusions follow, using propositions~\ref{Dpas0} and \ref{dualfenchel}.
	}
	\qed

The last statement of the previous proposition allows us to recover a result already stated
	in \cite{brasl76bcyi}.
\begin{coro}
	If for all $I \subset N_m$ of cardinality $n$, $d_I=\det((v_i))_{i \in I}$ is not zero,
	then for all $c=(c_i)_{i=1}^m$ such that:
	$$ \sum_{i=1}^m c_i =n \qquad \mathrm{and} \qquad 0 < c_i < 1  \mbox{ for all } i,$$
	then the \option{best}constant $D_c$ 
	\option{such that for all
	$m$-tuple of positive numbers $(\lambda_i)_{i=1}^m$, 
	$$ \det \left(\sum_{i=1}^m \lambda_i v_i \otimes v_i\right)
		\ge D_c \prod_{i=1}^m \lambda_i^{c_i} ,$$
	}is achieved for a certain $(\lambda_i)_{i=1}^m$.
\end{coro}
\option{
{\bf Proof:}
	It suffices to notice that because of the geometrical hypothesis
	on the vectors $(v_i)_{i=1}^m$, one has
	$$K=\mathrm{conv}\{1_I \, ; \; d_I \neq 0 \}
		=\left\{x \in [0,1]^m  \, ; \; \sum_{i=1}^m x_i=1 \right\}.$$
	The hypothesis on $c$ implies that it is in the
	relative interior of $K$.
	\qed
}

The following result shows that the reciprocal statement is almost true.
\begin{prop}
	If $(\R^n,(v_i)_{i=1}^m)$ is irreducible and if $c_1=1 $, then
	$D_c$ is achieved only when $m=n=1$. 
\end{prop}
\option{
{\bf Proof:}
	For $\lambda_i > 0$, $i=1, \ldots ,m$, one has
	\beq
		\Delta((\lambda_i)_{i=1}^m) &=& \frac{\sum_{|I|=n} \lambda_I d_I}
			{\prod_{i=1}^m \lambda_i^{c_i}} \\
		&=& \sum_{|I|=n, \, 1 \in I}  \frac{\lambda_{I-\{i\}} d_I}{\prod_{i=2}^m \lambda_i^{c_i}}
			+ \frac1{\la_1} \sum_{|I|=n, \, 1 \not\in I} \frac{\lambda_I d_I}
			{\prod_{i=2}^m \lambda_i^{c_i}} \cdot
	\eeq
	If  $\Delta$ achieves its minimum for positive $\la_i$, the coefficient
	of $\la_1$ has to be zero. This means that for all $I \subset N_m$ of 
	cardinality $n$ which does not contain 1, $(v_i)_{i \in I}$ is not a
	basis. It implies that rank$(v_2, \ldots, v_m)$ is less than $n-1$ and
	$$ \R^n = \R v_1 \oplus \mathrm{span}(v_2, \ldots, v_m).$$
	By irreducibility, it is possible only if $m=1$ and $n=1$.
	\qed
}

We come to unicity results: if $D$ is achieved, there is a unique minimizer, up to scalar 
	multiplication.
\begin{prop}
\label{uniqgauss}
	Assume that $(\R^n,(v_i)_{i=1}^m)$ has the irreducibility property. If $(\lambda_i)_{i=1}^m$
	and $(\mu_i)_{i=1}^m$ are two minimizers,\option{ i.e. if
	$$ \det \left(\sum_{i=1}^m \lambda_i v_i \otimes v_i\right)
		= D_c \prod_{i=1}^m \lambda_i^{c_i} ,$$
	and
	$$ \det \left(\sum_{i=1}^m \mu_i v_i \otimes v_i\right)
		=  D_c \prod_{i=1}^m \mu_i^{c_i} ,$$}
	then there exists $r \in \R$ such that for all $i$, $\lambda_i = r\mu_i$.
\end{prop}
{\bf Proof:}
	Let $t=((t_i)_{i=1}^m)$ and $s=((s_i)_{i=1}^m)$ such that for all $i$, one has
	$$ \lambda_i = e^{t_i} \mbox{  and  } \mu_i = e^{s_i}.$$
	Let $\psi$ be the function on $\R^m$ defined for all $((x_i)_{i=1}^m)$ by
	$$ \psi((x_i))= \phi ((x_i)) - \sum_{i=1}^m c_i x_i.$$
	Then $\psi $ reaches its minimum at the points $t$, $s$
	and also at $(t+s)/2$ because it is convex. So we have
	$$ \frac{ \phi(t)+\phi(s)}{2} = \phi \left( \frac{t+s}{2} \right),$$
	and there must be equality in the Cauchy-Schwartz inequality in the proof
	of proposition~\ref{lienDphi}. Hence, there exists $ a \in \R$ such that
	for all $I$, $|I|=n$,
	$$ d_I \exp\left( \sum_{i \in I} t_i \right) = a  \cdot d_I \exp\left( \sum_{i \in I} s_i \right).$$
	In particular, if $d_I \neq 0$, one has
	$$ \prod_{i \in I} \left( \frac{ \lambda_i}{\mu_i} \right) = a \, .$$
	Let $i,j \in N_m$ such that $i \bowtie j$; by definition, there exists $K \subset N_m$
	of cardinality $n-1$, such that $d_{\{i\} \cup K}$ and $d_{\{j\} \cup K}$ are both\
	non-zero. So, we have
	$$ \prod_{l \in \{i\} \cup K} \left( \frac{ \lambda_l}{\mu_l} \right)=
		\prod_{l \in \{j\} \cup K} \left( \frac{ \lambda_l}{\mu_l} \right), $$
	and after simplification
	$$ \frac{ \lambda_i}{\mu_i}= \frac{ \lambda_j}{\mu_j} \cdot$$
	By the irreducibility property (see Proposition~\ref{decomp}), this implies
	$\D \frac{ \lambda_1}{\mu_1}= \cdots =\frac{ \lambda_m}{\mu_m} \cdot $
	\qed

\subsection{The general case}
We have studied existence and uniqueness of centered Gaussian maximizers for (BL) and
	minimizers for (RBL), we turn to the general study. As explained before, we 
	may assume that $(\R^n, (v_i)_{i=1}^m)$ is irreducible. The behaviour of
	extremal functions is very different for $n=1$ and for $n \ge 2$.
\subsubsection{The case $n=1$}
If  $n=1$, then $n_i=1$ for all $i \le m$, the condition on $(c_i)_{i=1}^m$
	is just $\sum_{i=1}^m c_i=1$, and the $v_i$'s are just real numbers. The
	inequality (BL) is nothing else than H{\"o}lder's inequality for the 
	functions $x \mapsto f_i(v_i x)$, whereas (RBL) is the Pr\'ekopa-Leindler
	inequality for $x \mapsto f_i(x /v_i)$.

The equality cases can be settled from our proof; we will not do it because they are
	well-known: if $\sum_{i=1}^m c_i=1$, and $f_i \in L_1^+(\R)$, $i=1, \ldots, m$
	are non identically zero, then
	$$ \int_{\R} \prod_{i=1}^m f_i^{c_i} (x) \, dx =\prod_{i=1}^m \left( \int_{\R}  f_i
		\right)^{c_i}$$
	holds if and only if 
	$$ \frac{f_1}{\int_{\R}f_1}= \cdots =\frac{f_m}{\int_{\R}f_m} \cdot$$
	Under the same assumptions,
	$$  \int\limits_{\R}\sup\limits_{\sum c_i x_i =x} \prod_{i=1}^m f_i^{c_i} (x_i) \, dx 
		=\prod_{i=1}^m \left( \int_{\R}  f_i \right)^{c_i}$$
	holds if and only if there exists $(y_i)_{i=1}^m \in \R^m$ such that
	$$ \frac{f_1(\cdot - y_1)}{\int_{\R}f_1}= \cdots =\frac{f_m(\cdot - y_m)}{\int_{\R}f_m}
		\mbox{ is a log-concave function.}$$
 
\subsubsection{The case $n \ge 2$}
We prove that if there is a centered Gaussian extremizer, then up to dilatation and
	scalar multiplication, it is the only extremizer.
\begin{theo}
	Let $n \ge 2$ and let $(\R^n, (v_i)_{i=1}^m)$ be irreducible. Let $(c_i)_{i=1}^m$
	and $(\la_i)_{i=1}^m$ be positive numbers such that $D_c$ is achieved for 
	$(\la_i)_{i=1}^m$:
	$$ \det\left( \sum_{i=1}^m \lambda_i v_i \otimes v_i \right)= D_c \prod_{i=1}^m 
		\la_i^{c_i}.$$
	Then $(h_i)_{i=1}^m$ is a maximizer for (BL) if and only if there exist $a>0$, 
	$(\alpha_i)_{i=1}^m$
	positive and $ y \in \R^n$ such that for all $i$ and for all $t \in \R$,
	\begin{equation}
	\label{egbl}
		 h_i(t) = \alpha_i \exp ( - \la_i (at - \lag y, v_i \rag)^2) .
	\end{equation}
	The $m$-tuple $(h_i)_{i=1}^m$ is a minimizer for (RBL) if and only if there exist
	$b>0$, $(\beta_i)_{i=1}^m$ positive and $(t_i)_{i=1}^m$ real such that for all $i$ and for
	 all $t \in \R$,
	$$ h_i(t) = \beta_i \exp ( - (bt - t_i)^2/ \la_i) .$$
\end{theo}
{\bf Proof: }
	By Lemma~\ref{gaussbl} and by the proof of Lemma~\ref{gaussrbl}, we know that
	$(G_{\la_i})_{i=1}^m$ is a maximizer for (BL) and $(G_{\la_i^{-1}})_{i=1}^m$ a
	minimizer for (RBL), so by simple changes of
	variables in $\R^n$, one can check that the previous functions are
	extremizers.

	Let $(h_i)_{i=1}^m$ be a maximizer for (BL) and $(f_i)_{i=1}^m$ is a minimizer
	for (RBL). We may assume that $f_i,h_i$ are positive and continuous for all $i$. Indeed
	by the following lemma (which was communicated to me by K. Ball) we know that
	$(h_i \ast G_{\la_i})_{i=1}^m$ is a positive and continuous maximizer for (BL).
	If we know that it is Gaussian, then so is $(h_i)_{i=1}^m$ by the properties
	of the Fourier transform. The same argument is relevant for (RBL).
\begin{lem}

	If $(f_i)_{i=1}^m$ and $(g_i)_{i=1}^m$ are maximizers for (BL), then so is 
	$(f_i \ast g_i)_{i=1}^m$.

	If $(f_i)_{i=1}^m$ and $(g_i)_{i=1}^m$ are minimizers for (RBL), then so is 
	$(f_i \ast g_i)_{i=1}^m$.
\end{lem}
A proof of the first part of this lemma is written in \cite{bart98oyic}, the proof
	of the second part is very similar. Notice that this lemma is valid for
	the multidimensional version of the inequalities.

\vspace{3 mm}
We show now that if $(h_i)_{i=1}^m$ is a positive continuous maximizer for (BL), then it has
	to be of the form (\ref{egbl}); the proof for (RBL) is analogous and a bit simpler.
	We study precisely the proof of Lemma~\ref{fond} applied with $(h_i)_{i=1}^m$ being
	the maximizer we study and $(f_i)_{i=1}^m$ being the particular minimizer 
	for (RBL) that we know by hypothesis, namely
	$$ f_i(t) = \exp \left(- x^2 / \la_i \right) \, .$$
	Since our functions are positive, the change of variables $T_i$'s are 
	increasing differentiable bijections of $\R$, such that for all $t \in \R$,
	$$ T'_i (t)  \cdot f_i(T_i(t)) = h_i(t) .$$
	There must be equality in every step of the proof. In particular, for all $y \in \R^n$,
	one has
	$$ \det\left( \sum_{i=1}^m T'_i (\lag y, v_i \rag) v_i \otimes v_i \right)
		= D_c \prod_{i=1}^m  (T'_i (\lag y, v_i \rag))^{c_i}.$$
	By irreducibility and proposition~\ref{uniqgauss}, one gets for all $y \in \R^n$,
	$$ \frac{T'_1 (\lag y, v_i \rag)}{\la_1} = \cdots 
		= \frac{T'_m(\lag y, v_i \rag)}{\la_m} \cdot $$
	Since $n \ge 2$, for all $i \le m$ there exists $j \le m$ such that $v_i$  and
	$v_j$ are not colinear; so there exits $z \in \R$ such that
	$ \lag z, v_i \rag=1 \mbox{  and  } \lag z, v_i \rag=0.$
	The previous relation for $y=tz$ says that for all $t \in \R$
	$$ \frac{T'_i(t)}{\la_i} =\frac{T'_j(0)}{\la_j} \cdot$$
	Consequently, there exist $a>0$ and $(s_i)_{i=1}^m$ real such that for all $i$ and
	for all $t \in \R$,
	$$T'_i(t)=  a \la_i t + s_i \, $$
	and by the change of variable formula between $h_i$ and $f_i$ we get
	\beq
		h_i(t) &=& T'_i(t) \exp\left( -T_i^2(t)/\la_i \right) \option{ \\
		&=& a \la_i \exp\left(-(a \la_i t +s_i)^2 / \la_i \right) }  \\
		&=& \mu_i \exp( - \la_i(a t-t_i)^2 ) .
	\eeq
	for some positive $(\mu_i)$ and some real $(t_i)$.

	It remains to find which translates of a centered Gaussian maximizer 
	are still maximizers. Let $(g_i)_{i=1}^m$ be a maximizer,
	$$ g_i(t) = \exp( -\la_i t^2) ,$$
	and let $x=(x_i)_{i=1}^m \in \R^m$ and for $i\le m$, $h_i(t)= g_i(t-x_i)$.
	Let us consider $\R^m$ with the Euclidean metric given by
	$$ N^2(w)= \sum _{i=1}^m c_i \la_i w_i^2\, ,$$
	and the subspace
	$$ K=\left\{ (\lag y, v_i \rag)_{i=1}^m  \, ; \; y \in \R^n \right\} .$$
	Let $s$ be the orthogonal projection of $x$ onto $A$. Then 
	there exists $z\in \R^n$ satisfying $s_i=\lag z, v_i \rag$ for all $i$;
	moreover by the Pythagore Theorem
	$$ N((\lag y, v_i \rag-x_i)_{i=1}^m) \ge N((\lag y, v_i \rag-\lag z, v_i \rag)_{i=1}^m)$$
	with equality only if $x$ belongs to $A$, that is $x=s$.
	\option[Thus $J((h_i)_{i=1}^m) \le J((g_i)_{i=1}^m)$, with equality]{
	Since
	\beq
		J((h_i)_{i=1}^m)&=& \int_{\R^n} \exp\left(- \sum_{i=1}^m c_i \la_i
			(\lag y, v_i \rag-x_i)^2 \right) \, dy \\
		&=& \int_{\R^n} \exp\left(- N^2((\lag y, v_i \rag-x_i)_{i=1}^m) \right)\, dy \\
		& \le& \int_{\R^n} \exp\left(- N^2((\lag y, v_i \rag- \lag z, v_i \rag)_{i=1}^m)
				 \right)\, dy \\
		&=& \int_{\R^n} \exp\left(- \sum_{i=1}^m c_i \la_i
			(\lag y-z, v_i \rag-x_i)^2 \right) \, dy \\
		&=& J((g_i)_{i=1}^m)
	\eeq
	with equality only if $x \in A$, $(h_i)_{i=1}^m$ is a maximizer
	} 
	only if $x_i=\lag z, v_i \rag$
	for all $i$. 
	\qed

\vspace{4 mm}
{\bf Remark: }
There are, for $n \ge 2$, some remaining questions. If there is no centered Gaussian maximizer,
	is there any maximizer at all? The answer seems to be no: if $(f_i)_{i=1}^m$ is a maximizer
	then by the Brascamp-Lieb-Luttinger inequality \cite{brasll74grim} so is $(f^*_i)_{i=1}^m$,
	where $f^*$ is the symmetric rearrangement. As K. Ball remarked, for all $k$ integer
	$$ \Big( \sqrt k \underbrace{f_i^{*} \ast \cdots \ast f_i^*}_{k \mbox{ times}}(\sqrt k \,  \cdot )
		\Big)_{i=1}^m $$
	is also a maximizer; moreover, but under some integrability assumptions, it tends towards
	a centered Gaussian $m$-tuple by the Central Limit Theorem.

Notice that our method gives the answer when there are positive continuous maximizers
	for (BL) and positive continuous minimizers for the corresponding (RBL). The
	study the equality case of Lemma~\ref{fond} shows that the constant $D$ must be achieved,
	so there is a centered Gaussian maximizer.

Some arguments based on the equality case in the Minkowski inequality (see \cite{lieb90gkgm})
	might help to solve the question.

\section{Applications to convex geometry}
\subsection{Dimension one}
	K. Ball noticed that an additionnal geometrical hypothesis on the vectors $(v_i)$,
	which is frequent in convexity, allows an easy computation of the optimal constants
	in the Brascamp-Lieb inequality. For completeness, we begin by the proof of his
	observation.
\begin{prop}
\label{mingaussball}
	Let $m\ge n$, let $v_1, \ldots, v_m$ be vectors in $\R^n$ such that
	$ \sum_{i=1}^m v_i \otimes v_i =I_n,$
	where $I_n$ stands for the identity map; then for every $m$-tuple $(\lambda_i)_{i=1}^m$
	of positive numbers
	$$  \det \left (\sum_{i=1}^m \lambda_i v_i \otimes v_i \right) \ge \prod_{i=1}^m 
		\lambda_i^{|v_i|^2}.$$
	There is equality when $\lambda_1= \cdots = \lambda_m$ .
\end{prop}
{\bf Proof:}
	By the Cauchy-Binet formula, we have
	$$ 1= \det I_n = \det\left( \sum_{i=1}^m v_i \otimes v_i \right)
		= \sum_{|I|=n}  d_I. $$
	Hence we can use the arithmetic-geometric inequality with coefficients $d_I$ :
	$$\det\left( \sum_{i=1}^m \lambda_i v_i \otimes v_i \right)
		= \sum_{|I|=n} \lambda_I d_I \ge \prod_{|I|=n} \lambda_I^{d_I}.$$
	Each $\lambda_i$ appears with the total exponent $\sum_{I, \, i \in I} d_I$.
	But by Corollary~\ref{cb} applied to the $m$-tuple 
	$(v_1,\ldots, v_{i-1}, 0, v_{i+1}, \ldots , v_m)$ we get:
	\beq
		\sum_{I, \, i \in I} d_I &=& \sum_{I} d_I - \sum_{I, \, i \not\in I} d_I \\
		&=& 1- \det (v_1 \otimes v_1 + \cdots+ v_{i-1}\otimes v_{i-1} +v_{i+1}\otimes v_{i+1}
			+ \cdots + v_m\otimes v_m )  \\
		&=& 1- \det (I_n -v_i\otimes v_i) = |v_i|^2.
	\eeq
	\option{
	So we have proved 
	$$  \det \left (\sum_{i=1}^m \lambda_i v_i \otimes v_i \right) \ge \prod_{i=1}^m 
		\lambda_i^{|v_i|^2}.$$
	When all the $\lambda_i$'s are equal this becomes an equality because $\sum_{i=1}^m |v_i|^2=n$
	(take the traces in $ \sum_{i=1}^m v_i \otimes v_i =I_n$).
	}
	\qed

\vspace{4 mm}
\option{
The additional hypothesis also gives more structure to the adapted decomposition:
\begin{lem}
	Let $v_1, \ldots, v_m$ in $\R^n$, let
	$C$ be the best adapted partition  built in proposition~\ref{decomp}. If
	$ \sum_{i=1}^m v_i \otimes v_i =I_n$, then $C$ induces an orthonormal
	decomposition of the space, namely
	$$ \R^n = \bigoplus_{I \in C}^{\bot} E_I \, ,$$
	where $E_I= \mathrm{ span }\{ (v_i)_{i \in I} \}$.
\end{lem}
{\bf Proof:}
	By proposition~\ref{decomp}, we just know 
	$$\R^n = \bigoplus_{I \in C} E_I \, ,$$
	so we have a family $(P_I)_{I \in C}$ of projections such that for all $I$,
	$\mathrm{Im}(P_I)=E_I$ and $\ker (P_I)=\bigoplus_{J\neq I} E_J$.
	From the decomposition of the identity $ \sum_{i=1}^m v_i \otimes v_i =I_n$,
	we obtain that for all $x$ in $\R^n$,
	$$ x= \sum_{i=1}^n \langle x, v_i \rangle v_i = \sum_{I \in C} \left(
		\sum_{i \in I} \langle x, v_i \rangle v_i \right). $$
	By unicity of the decomposition of $x$ in the direct sum induced by $C$, we 
	have for all $I \in C$:
	$$P_I (x) =	\sum_{i \in I} \langle x, v_i \rangle v_i \, ,$$
	so $P_I$ is an orthogonal projection. It is clear that the direct sum $ \bigoplus_{I \in C} E_I $
	must be  orthogonal.
	\qed
}

Ball's version of the Brascamp-Lieb inequality and the corresponding reverse version are 
	as follows:
\begin{theo}
\label{blball}
	Let $m \ge n$, let $(u_i)_{i=1}^m$ be unit vectors in $\R^n$ and let $(c_i)_{i=1}^m$
	be positive real numbers such that
	$$ \sum_{i=1}^m c_i u_i \otimes u_i = I_n .$$
	Then for all $f_i \in L_1^+( \R)$, $i=1, \ldots, m$ one has
	$$ \int_{\R^n} \prod_{i=1}^m f_i^{c_i}( \langle x, u_i \rangle) \, dx 
		\le \prod_{i=1}^m \left( \int f_i \right) ^{c_i},$$
	and 
	$$  \int\limits_{\R^n}^*  \sup_{x=\sum c_{i} \theta_i u_{i}} 
		\prod_{i=1}^{m} f_{i}^{c_i}(\theta_{i}) \, dx
		\ge \prod_{i=1}^m \left( \int f_i \right)^{c_i}.$$
\end{theo} 
\option{
This form of the Brascamp-Lieb inequality has several applications in convex geometry, so
	does the reverse form. We give here as a corollary a lower estimate for the volume
	of zonoids, already proved by Ball \cite{ball91scb} by a geometrical inductive method.
\begin{coro}
	Let $m \ge n$, let $(u_i)_{i=1}^m$ be unit vectors in $\R^n$ and let $(c_i)_{i=1}^m$
	be positive real numbers such that
	$$ \sum_{i=1}^m c_i u_i \otimes u_i = I_n \, .$$
	Then for every $\alpha_1, \ldots, \alpha_m$ positive,
	$$ \mathrm{Vol}\left( \sum_{i=1}^m \alpha_i [-u_i,u_i] \right)
		\ge 2^n \prod_{i=1}^m \left( \frac{\alpha_i}{c_i} \right)^{c_i} .$$
\end{coro}
{\bf Proof:}
	Apply the reverse form when $f_i$ is the characteritic function of the interval
	$]- \alpha_i / c_i, +\alpha_i / c_i [$.
	\qed

For more sophisticated applications we refer to \cite{bart98}.
}

The equality case is completly settled: the space $\R^n$ is an orthogonal sum of
	irreducible subspaces. On the irreducible subspaces of dimension one there
	is equality for (BL) if and only if the functions are equal up to scalar multiplication
	and, for (RBL), if and only if all the functions are equal up to multiplication and 
	translation to a common log-concave function. On irreducible spaces of dimension
	more than or equal to 2, there is equality if and only if the functions are
	(up to scalar multiplications, up to translations and only coherent translations
	in the direct form) equal to a common centered Gaussian function. We state a usefull corollary
	which allows to solve the equality case in the geometric applications due to K. Ball.
\begin{coro}
	Let $m \ge n$, let $(u_i)_{i=1}^m$ be $m$ different unit vectors in $\R^n$ and
	 let $(c_i)_{i=1}^m$ be positive real numbers such that
	$$ \sum_{i=1}^m c_i u_i \otimes u_i = I_n \, .$$
	If $(f_i)_{i=1}^m$ are non-identically-zero functions in $L_1^+(\R)$ such that none of them
	is a Gaussian and
	$$\int_{\R^n} \prod_{i=1}^m f_i^{c_i}( \langle x, u_i \rangle) \, dx 
		= \prod_{i=1}^m \left( \int f_i \right) ^{c_i},$$
	then $(u_i)_{i=1}^m$ is an orthonormal basis of $\R^n$.
\end{coro} 
We give an application: for $K \subset \R^n$ a convex body, we denote by $\mathrm{vr}(K)$ the
	ratio of the volume of $K$ by the volume of the maximal volume ellipsoid contained in $K$
	(called the John's ellipsoid, see \cite{john48epis}). In \cite{ball89vscr} and
	\cite{ball91vrri}, K. Ball proved by means of the Brascamp-Lieb inequality that
	simplices have maximal volume ratio and that among symmetric bodies, parallelotopes have.
	By the previous corollary, we can answer the question of equality cases. We denote by $Q_n$
	the unit cube and by $\Delta_n$ the regular simplex.
\begin{prop}
	Let $K \subset \R^n$ a convex body. 
	\begin{itemize}
		\item If $K$ is symmetric and $\mathrm{vr}(K)=\mathrm{vr}(Q_n) $ then $K$ is a parallelotope.
		\item If $\mathrm{vr}(K)=\mathrm{vr}(\Delta_n) $ then $K$ is a simplex.
	\end{itemize}
\end{prop}

\subsection{Larger dimensions}
We obtain a multidimensonal generalization of Ball's version of the Brascamp-Lieb inequality and
	its converse. The estimate for Gaussians is a generalization of proposition~\ref{mingaussball}.
\option{
\begin{prop} 
	For $i\le m$ let $B_i \in \mathcal L (\R^n,\R^{n_i})$ such that for all $x \in \R^{n_i}$,
	$ |B_i^*(x)|=|x|$. Assume that we have some positive numbers $(c_i)_{i=1}^m$ satisfying
	$$ \sum_{i=1}^m c_i B_i^*B_i = I_n $$
	Then for all $A_i \in \mathcal S ^{+}(\R^{n_i})$,
	$1=1, \ldots, m$, one has
	$$ \det \left( \sum_{i=1}^m c_i B_i^* A_i B_i \right) \ge \prod_{i=1}^m
		\left( \det A_i \right)^{c_i} .$$
\end{prop}
{\bf Proof:}
	For each $i$, $A_i$ can be written as
	$$ A_i = \sum_{j=1}^{n_i} \lambda_{i,j} e_{i,j} \otimes e_{i,j} \, ,$$
	where the vectors $(e_{i,j})_{j=1}^{n_i}$ are an orthonormal basis of 
	$\R^{n_i}$ formed of eigenvectors of $A_i$ and $(\lambda_{i,j})_{j=1}^{n_i}$
	are its eigenvalues.
	Since for all $i$, 
	$$ I_{n_i}= \sum_{j=1}^{n_i} e_{i,j} \otimes e_{i,j} \, ,$$
	one has
	$$ \sum_{i,j} c_i B_i^* e_{i,j} \otimes B_i^* e_{i,j}= I_n \, .$$
	and the vectors $B_i^* e_{i,j}$ have norm one. Thus we can apply the result of
	proposition~\ref{mingaussball} and obtain a lower estimate for the determinant
	of
	$$ \sum_{i=1}^m c_i B_i^* A_i B_i = \sum_{i=1}^m \sum_{j=1}^{n_i} c_i \lambda_{i,j}
		B_i^* e_{i,j} \otimes B_i^* e_{i,j} \; .$$
	The determinant is bounded below by 
	$$ \prod_{i=1}^m \prod_{j=1}^{n_i} \lambda_{i,j} ^{c_i}
		 = \prod_{i=1}^m \left(\det A_i \right)^{c_i} .$$
	\qed

As a consequence of the main theorem and of the previous estimate, one has
}
\begin{theo}
\label{multballBL}
	Let $m,n$ be integers. For $i=1, \ldots, m$ let $E_i$ be a subspace of $\R^n$ of 
	dimension $n_i$ and let $P_i$ be the orthogonal projection onto $E_i$ (on each $E_i$ there
	is a Lebesgue measure compatible with the induced Euclidean structure). Assume that
	there exist positive numbers $(c_i)_{i=1}^m$ satisfying
	$$ \sum_{i=1}^m c_i P_i =I_n \, ,$$
	Then if for $i=1, \ldots, m$, $f_i$ is a non-negative  integrable function on $E_i$, one has
	$$ \int_{\R^n} \prod_{i=1}^m f_i^{c_i}(P_i x) \, d^n x \le
		 \prod_{i=1}^m \left(\int_{E_i} f_i \right)^{c_i}, $$
	and 
	$$ \int\limits_{\R^n}^* \sup\limits_{x=\sum_{i=1}^m c_i x_i, \, x_i \in E_i} 
		\prod_{i=1}^m f_i^{c_i}(x_i) \, d^n x 
		\ge \prod_{i=1}^m \left(\int_{E_i} f_i \right)^{c_i}. $$
\end{theo}
{\bf Remark:} When the $f_i$ are taken to be characteristic functions of sets, the reverse
	inequality provides a Brunn-Minkowski type result for convex bodies which do not
	have full dimension: if $K_i \subset E_i$ then
	$$ \mathrm{Vol}_n\left( \sum_{i=1}^m c_i K_I \right) 
		\ge \prod_{i=1}^m \left(\mathrm{Vol}_{E_i}(K_i)\right)^{c_i}.$$
	\option[When the sets $K_i$ are segments, we recover a lower estimate for the volume
	of zonoids, already proved by Ball \cite{ball91scb} by a geometrical inductive method.]{}
\section{Link with Young's convolution inequality}
In \cite{bart98oyic}, the author gave a proof of Young's inequality and its converse
	which is also based on measure-preserving mappings. Following an idea of 
	\cite{brasl76bcyi}, we illustrate how a generalization of Young's inequality
	and its converse contains Ball's version of (BL) and the corresponding form
	of (RBL). For $t>1$, we define the conjugate 
	number $t'$ by $1/t + 1/t' =1$.
\begin{theo}
\label{geneyoung}
	Let $m \ge n$ be integers, let $V$ be an orthogonal $m \times m$ matrix and denote
	by $(v_i)_{i=1}^m$ its rows. Let $M$ be the $(m-n) \times m$ submatrix of $V$
	formed of the last $m-n$ columns. Let $r$ and $(p_i)_{i=1}^m$ be larger than 1 and 
	such that $\sum_{i=1}^m 1/p_i = n + (m-n)/r$;  let
	$D=D_{r,p_i}$ be the largest constant such that for all positive $(\la_i)_{i=1}^m$, one has
	$$ \det \left( M \cdot \mathrm{diag}(\la_i) \cdot ^t \!M \right) \ge D \prod_{i=1}^m 
		\la_i^{r'/p'} .$$
	Then for every continuous positive integrable functions on $\R$, $(f_i)_{i=1}^m$ and
	$(F_i)_{i=1}^m$ such that for all $i$, $\int f_i = \int F_i$, one has
	$$ \left[ \int \left[ \int  \prod_{i=1}^m f_i^{1/p_i}( \lag x, v_i \rag )
		\, dx_1 \cdots dx_n \right]^r  dx_{n+1} \cdots dx_m \right]^{1/r}  \qquad \qquad$$
	$$ \qquad \qquad \le D^{-1/r'} \int \left[ \int \prod_{i=1}^m F_i^{r/p_i}( \lag X, v_i \rag )
		\, dX_m \cdots dX_{n+1} \right]^{1/r}  dX_n \cdots dX_1 \, . $$
\end{theo}
{\bf Proof:} For all $i$, there exists a positive differentiable increasing map $T_i$ 
	satisfying, for all $s \in \R$
	$$ \int_{- \infty}^{T_i(s)} f_i = \int_{- \infty}^s F_i \, $$
	and by differentiation:
	\begin{equation}
	\label{reldiff}
		T'_i(s) f_i(T_i(s)) = F_i(s) \, .
	\end{equation}
	We consider the change of variable $\Theta$ in $\R^m$ given by $V^{-1} (T_1 \otimes \cdots
	\otimes T_m) V$. More precisely, $x= \Theta (X)$ means that for all $i$,
	$$ \lag v_i , x \rag = T_i (\lag v_i , X \rag ). $$
	The application $\Theta$ is clearly bijective, its differential at a point $X \in \R^n$
	is 
	$$ d\Theta(X)= V^{-1} \mathrm{diag}( T'_i( \lag v_i , X \rag )) V = ^t \! 
		V \mathrm{diag}( T'_i( \lag v_i , X \rag )) V \, , $$
	so its jacobian is simply
	$$ \prod_{i=1}^m T'_i( \lag v_i , X \rag ) \, .$$
	We want an upper estimate of the integral
	$$ I=  \left[ \int \left[ \int  \prod_{i=1}^m f_i^{1/p_i}( \lag x, v_i \rag )
		\, dx_1 \cdots dx_n \right]^r  dx_{n+1} \cdots dx_m \right]^{1/r} \, ,$$
	which is finite (we may suppose that all our functions are dominated by 
	some Gaussian function). Hence there exists a positive function $h \in L^{r'}(\R^{m-n})$
	such that $\|h\|_{r'}=1$ and
	$$ I=  \int  \prod_{i=1}^m f_i^{1/p_i}( \lag x, v_i \rag ) h(x_{n+1}, \ldots, x_m)
		\, dx_1 \cdots dx_m  \, .$$
	By the change of variables $x= \Theta (X)$ and by relations (\ref{reldiff}), we get
	$$ I =  \int  \prod_{i=1}^m f_i^{1/p_i}(T_i( \lag x, v_i \rag )) h(A(x)) 
		\prod_{i=1}^m T'_i( \lag v_i , X \rag ) \, dX_1 \cdots dX_m$$
	$$ = \int \left[ \int \prod_{i=1}^m F_i^{1/p_i}( \lag X, v_i \rag )h(A(x))
		\prod_{i=1}^m \big(T'_i( \lag v_i , X \rag )\big)^{1/p'_i} \, dX_m \cdots dX_{n+1}
		\right]  dX_n \cdots dX_1 \, , $$
	where $A$ is defined by
	$$ A(x)= \left( \sum_{i=1}^m M_{n+1,i} T_i(\lag v_i , X \rag ), \ldots ,
			 \sum_{i=1}^m M_{m,i} T_i(\lag v_i , X \rag )  \right)\, .$$
	For fixed $X_1, \ldots, X_n$, we consider $A$ as a function of $X_{n+1}, \ldots, X_m$.
	Its differential is 
	$$ M \cdot \mathrm{diag} (T'_i( \lag v_i , X \rag  )) \cdot ^t \!\! M \, ,$$
	so we know, by hypothesis, a lower estimate for its Jacobian
	$$ \det dA(X) \ge D \prod_{i=1}^m \bigg(T'_i( \lag v_i , X \rag )\bigg)^{r'/p'_i} \, .$$
	By H{\"o}lder's inequality with parameters $(r,r')$ applied to the inner integral
	of the previous expression for $I$ and by the lower estimate for the Jacobian of $A$,
	we get
	$$ I \le D^{-1/r'} \int \left[  \int \prod_{i=1}^m F_i^{r/p_i}( \lag X, v_i \rag )
		  \, dX_m \cdots dX_{n+1}  \right]^{1/r} \hspace{5 cm} $$
	$$ \hspace{5 cm} \left[ \int h^{r'}(A(X)) \det(dA(X)) 
		  \, dX_m \cdots dX_{n+1} \right]^{1/r'} dX_n \cdots dX_1 \,.$$
	Since for fixed $X_1, \ldots, X_n$, $A$ is injective (indeed its differential is
	symmetric definite positive),
	$$ \int h^{r'}(A(X)) \det(dA(X)) \, dX_m \cdots dX_{n+1} \le \int h^{r'} =1 \, ,$$
	so
	$$ I \le  D^{-1/r'} \int \left[  \int \prod_{i=1}^m F_i^{r/p_i}( \lag X, v_i \rag )
		  \, dX_m \cdots dX_{n+1}  \right]^{1/r} \, dX_n \cdots dX_1 \,.$$
	\qed

This theorem contains both (BL) and (RBL) in the form stated in Theorem~\ref{blball}. We have 
	positive $(c_i)_{i=1}^m$ and unit vectors $(u_i)_{i=1}^m$ in $\R^n$, linked by
	$$ \sum_{i=1}^m c_i u_i \otimes u_i =I_n\, .$$
	We begin by a very simple fact:
\begin{lem}
	Let $m \ge n$ and let $(c_i)_{i=1}^m$ be positive $(u_i)_{i=1}^m$ be unit vectors in $\R^n$
	with the relation:
	$$ \sum_{i=1}^m c_i u_i \otimes u_i =I_n\, ,$$
	then there exists an orthonormal basis $(v_i)_{i=1}^m$ of $\R^m$ such that for all $i$,
	$$ P(v_i) = \sqrt{c_i} u_i \, ,$$
	where $P$ stands for the projection from $\R^m$ onto $\R^n$ which keeps only the first
	$n$ coordinates of a vector.
\end{lem}
We apply Theorem~\ref{geneyoung} with the vectors $(v_i)_{i=1}^m$ given by the lemma and for 
	special values of the parameters: namely for $R>1$ and very close to one, we chose
	$$ p_i =\frac1{Rc_i} $$
	and 
	$$ r= \frac{m-n}{\sum_{i=1}^m 1/p_i - n}= \frac{m-n}{n(R-1)}.$$
	Then we take the limits when $R$ tends to $1$ in the inequality provided by the 
	theorem.

	Let us describe the asymptotic behaviour of all the related quantities when
	$R \to 1$. It is clear that
	$$ \frac1{p_i}  \to c_i \, ,  \qquad \frac1{p'_i}  \to 1-c_i \, ,  \qquad r \to \infty 
		\qquad \mathrm{and} \qquad r' \to 1.$$
	We define $m$ vectors $(w_i)_{i=1}^m$ in $\R^{m-n}$ as follows: the coordinates of
	$w_i$ are the last coordinates of $v_i$. With these notations, $D_{r,p_i}$
	is the largest constant such that for all positive $(\la_i)_{i=1}^m$,
	$$ \det \left( \sum_{i=1}^m \la_i w_i \otimes w_i \right) \ge D_{r,p_i} 
		\prod_{i=1}^m \la_i ^{r'/p'_i} .$$
	Since we have the orthogonal decomposition $v_i = \sqrt{c_i} u_i + w_i$, we know
	that $|w_i|^2=1-c_i$. Moreover, $(w_i)_{i=1}^m$ being the orthogonal projection of
	and orthonormal basis, the following relation holds
	$$ \sum_{i=1}^m  w_i \otimes w_i =I_{m-n} \, ,$$
	so we get from Lemma~\ref{mingaussball}, for all positive $(\la_i)_{i=1}^m$,
	$$ \det \left( \sum_{i=1}^m \la_i w_i \otimes w_i \right) \ge 
		\prod_{i=1}^m \la_i ^{1-c_i} .$$
	As $r'/p'_i \to 1-c_i$ when $R$ tends to one it follows that $D_{r,p_i} \to 1$.
	
We study now the quantities involving the functions. The main point is that when a function 
	decreases fast enough, its $r$-norm tends to its essential supremum when $r$ tends
	to infinity. We introduce some more notation: for $(x_1, \ldots, x_m)$ we set 
	$y=(x_1, \ldots, x_n)$ and $z=(x_{m+1}, \ldots, x_m)$.
	When $R \to 1$, $r \to \infty$, and
	$$  \left[ \int \left[ \int  \prod_{i=1}^m f_i^{1/p_i}( \lag x, v_i \rag )
		\, dx_1 \cdots dx_n \right]^r  dx_{n+1} \cdots dx_m \right]^{1/r} $$
	tends to 
	$$ \sup\limits_{z \in \R^{m-n}} \int_{\R^n} \prod_{i=1}^m f_i^{c_i}
		( \sqrt{c_i} \lag y, u_i \rag + \lag z, w_i \rag) \, dy $$
	which is larger than
	$$ \int_{\R^n} \prod_{i=1}^m f_i^{c_i} ( \sqrt{c_i} \lag y, u_i \rag) \, dy \, .$$
	On the other hand, the quantity 
	$$\int \left[ \int \left( \prod_{i=1}^m F_i^{1/p_i}( \lag x, v_i \rag )\right)^r
		 dx_m \cdots dx_{n+1} \right]^{1/r} \!  dx_n \cdots dx_1 $$
	tends to 
	$$ \int  \sup\limits_{z}\prod_{i=1}^m F_i^{c_i}( \lag x, v_i \rag)\,  dx_n \cdots dx_1 \, ,$$
	where the supremum is on the vectors $z=(x_{n+1}, \ldots, x_m) \in \R^{m-n}$. 
	Noticing that the numbers $\lag x, v_i \rag$ are just the coordinates of $x=y+z$ in
	the orthonormal basis $(v_i)_{i=1}^m$, and since the existence of a $z \in \R^{m-n}$
	such that $y+z= \sum_{i=1}^m \alpha_i v_i$ is equivalent to $y=\sum_{i=1}^m \sqrt{c_i}
	\alpha_i u_i$ (by taking orthogonal projection onto the first $n$ coordinates),
	the previous integral is just
	$$ \int\limits_{\R^n}  \sup\limits_{\sum \sqrt{c_i} \alpha_i u_i =y}
		\prod_{i=1}^m F_i^{c_i}( \alpha_i) \, dy \, ,$$
	$$ =\int\limits_{\R^n}  \sup\limits_{\sum c_i \theta_i u_i =y}
		\prod_{i=1}^m F_i^{c_i}( \sqrt{c_i} \theta_i) \, dy \, .$$
	Hence the limiting case of Theorem~\ref{geneyoung} states that for all $f_i, F_i$, one has
	$$ \frac{ \D \int\limits_{\R^n} \prod_{i=1}^m f_i^{c_i} (  \lag y, u_i \rag) \, dy}
		{\D \prod_{i=1}^m  \left( \int f_i\right) ^{c_i}  }
		\le
	   \frac{\D	\int\limits_{\R^n}  \sup\limits_{\sum c_i \theta_i u_i =y}
		\prod_{i=1}^m F_i^{c_i}( \theta_i)  \, dy}
		{\D \prod_{i=1}^m  \left( \int F_i\right) ^{c_i}} \, \cdot$$
	When the $F_i$'s are identical centered Gaussians this is Ball's version of (BL),
	and when the $f_i$'s are identical centered Gaussians it is the corresponding 
	version of (RBL). 
	\qed
	

\signature

\end{document}